\documentclass[12pt]{article}
\usepackage{latexsym}
\usepackage{amssymb}
\usepackage{graphicx}
\usepackage{cite}
\usepackage{color}
\usepackage{ifpdf}

\newtheorem{Definition}{Definition}[part]
\newtheorem{Proposition}{Proposition}[part]

\newtheorem{Lemma}{Lemma}[part]
\newtheorem{Corollary}{Corollary}[part]

\newtheorem{Conjecture}{Conjecture}

\def \ep{\hbox{ }\hfill$\Box$}

\def\reff#1{{\rm(\ref{#1})}}

\begin{document}
\title{Cored Hypergraphs, Power Hypergraphs and Their Laplacian H-Eigenvalues\thanks{This
work was supported by the Hong Kong Research Grant
Council (Grant No. PolyU 501909, 502510, 502111 and 501212) and NSF of China (Grant No. 11231004).}}

\author{
Shenglong Hu \thanks{Email: Tim.Hu@connect.polyu.hk. Department of
Applied Mathematics, The Hong Kong Polytechnic University, Hung Hom,
Kowloon, Hong Kong.},\hspace{4mm} Liqun Qi \thanks{Email:
maqilq@polyu.edu.hk. Department of Applied Mathematics, The Hong
Kong Polytechnic University, Hung Hom, Kowloon, Hong Kong.},
\hspace{4mm} Jia-Yu Shao \thanks{Email: jyshao@tongji.edu.cn. Department of Mathematics, Tongji University, Shanghai, China.}}

\date{\today}
\maketitle

\begin{abstract}
In this paper, we introduce the class of cored hypergraphs and power
hypergraphs, and investigate the properties of their Laplacian
H-eigenvalues.   From an ordinary graph, one may generate a
$k$-uniform hypergraph, called the $k$th power hypergraph of that
graph.   Power hypergraphs are cored hypergraphs, but not vice
versa. Hyperstars, hypercycles, hyperpaths are special cases of
power hypergraphs, while sunflowers are a subclass of cored
hypergraphs, but not power graphs in general. We show that the
largest Laplacian H-eigenvalue of an even-uniform cored hypergraph
is equal to its largest signless Laplacian H-eigenvalue. Especially,
we find out these largest H-eigenvalues for even-uniform sunflowers.
Moreover, we show that the largest Laplacian H-eigenvalue of an
odd-uniform sunflower, hypercycle and hyperpath is equal to the
maximum degree, i.e., $2$.  We also compute out the H-spectra of
the class of hyperstars. When $k$ is odd, the H-spectra of the hypercycle of size $3$ and the
hyperpath of length $3$ are characterized as well.

\vspace{2mm}
\noindent {\bf Key words:}\hspace{2mm} Tensor, H-eigenvalue, hypergraph, Laplacian, power hypergraph \vspace{1mm}

\noindent {\bf MSC (2010):}\hspace{2mm}
05C65; 15A18
\end{abstract}

\section{Introduction}
\setcounter{Theorem}{0} \setcounter{Proposition}{0}
\setcounter{Corollary}{0} \setcounter{Lemma}{0}
\setcounter{Definition}{0} \setcounter{Remark}{0}
\setcounter{Conjecture}{0}  \setcounter{Example}{0} \hspace{4mm}
A natural definition for the Laplacian tensor and the signless
Laplacian tensor of a $k$-uniform hypergraph for $k \ge 3$ was
introduced in \cite{q12a}.  See Definition 2.2 of this paper.
Recently, Hu, Qi and Xie \cite{hqx13} studied the largest Laplacian
and signless Laplacian eigenvalues of a $k$-uniform hypergraph, and
generalized some classical results of spectral graph theory to
spectral hypergraph theory, in particular, when $k$ is even.

One classical result in spectral graph theory \cite{bh11, z07} is
that the largest Laplacian eigenvalue of a graph is always less than
or equal to the largest signless Laplacian eigenvalue, and when the
graph is connected, the equality holds if and only if the graph is
bipartite.

In \cite{hqx13}, it was shown that the largest Laplacian
H-eigenvalue of a $k$-uniform hypergraph is always less than or
equal to the largest signless Laplacian H-eigenvalue, and when the
hypergraph is connected, the equality holds if and only if the hypergraph is
odd-bipartite.  A $k$-uniform hypergraph is called odd-bipartite if
$k$ is even, and the vertex set of the hypergraph can be divided to two
parts, such that each edge has odd number of vertices in each of
these two parts.

Hu, Qi and Xie \cite{hqx13}  generalized cycles in graphs to
hypercycles in $k$-uniform hypergraphs.  By showing that the largest
signless Laplacian H-eigenvalue for a hypercycle is computable and
an even-order hypercycle is odd-bipartite, they showed that the
largest Laplacian H-eigenvalue of a hypercycle is computable when
$k$ is even.

Another classical result in spectral graph theory \cite{bh11, gm94,
z07} is that the largest Laplacian eigenvalue of a graph is always
greater than or equal to the maximum degree of that graph, plus one.
The lower bound is attained if there exists a vertex adjacent to all
the other vertices of that graph.

Hu, Qi and Xie \cite{hqx13} showed that when $k$ is even the largest
Laplacian H-eigenvalue of a $k$-uniform hypergraph is always greater
than or equal to the maximum degree of that hypergraph, plus $\alpha(k)$,
where $\alpha(k) > 0$ and $\alpha(2) = 1$. They showed that the
lower bound is attained if the hypergraph is a hyperstar.

When $k$ is odd, the situation is very different.   It was shown in
\cite{hqx13} that in this case the largest Laplacian H-eigenvalue is
always strictly less than the largest signless Laplacian
H-eigenvalue, and the lower bound of the largest Laplacian
H-eigenvalue is the maximum degree itself, attained by any
hyperstar.

The results of \cite{hqx13} raised the interests to study the
Laplacian H-eigenvalues of a $k$-uniform hypergraph.   Actually,
they posed several questions for further research.   First, when $k$ is odd, is a hyperstar the only example of a $k$-uniform
hypergraph whose largest Laplacian H-eigenvalue is equal to the
maximum degree?   Second, can we calculate all Laplacian
H-eigenvalues for some special $k$-uniform hypergraphs, such as
hyperstars and hypercycles?   This is useful if one wishes to study
the second smallest Laplacian H-eigenvalue (with multiplicity) of a $k$-uniform
hypergraph, as the second smallest Laplacian eigenvalue of a graph
plays a key role in spectral graph theory \cite{bh11, c97, z07}.

Motivated by these questions, we study Laplacian H-eigenvalues of some
special $k$-uniform hypergraphs in this paper.

With the same way to generalize stars and cycles to hyperstars to
hypercycles, we may generalize an arbitrary graph $G$ to a
$k$-uniform hypergraphs $G^k$.  See Definition 2.4 of this paper. We
call $G^k$ the $k$th power of $G$, hence call it a power hypergraph.
In particular, paths are generalized to hyperpaths.   We will see
that when $k$ is even, a power hypergraph is odd-bipartite.   We may
conclude this for a broader class of $k$-uniform hypergraphs. We
call such hypergraphs cored hypergraphs.  See Definition 2.3 of this
paper. A power hypergraph is a cored hyoergraph but not vice versa.
In particular, we introduce a special subclass of cored hypergraphs,
called sunflowers.   See Definition 3.1 of this paper. A sunflower
is not a power hypergraph in general.   We show that when $k$ is
even, a cored hypergraph is odd-bipartite.   Thus, when $k$ is even,
the largest Laplacian H-eigenvalue and the largest signless Laplcian
H-eigenvalue of a cored hypergraph is the same. This enhances our
understanding on odd-bipartite hypergraphs and their largest
Laplacian eigenvalues.   We will show that the largest Laplacian
H-eigenvalue of an even-order sunflower is computable.

Then, when $k$ is odd, we will show that the largest Laplacian
H-eigenvalue of an odd-uniform sunflower, hypercycle and hyperpath
is equal to the maximum degree, i.e., $2$.   This shows that for a
very broad class of hypergraphs, when $k$ is odd, the largest
Laplacian H-eigenvalue is equal to the maximum degree of the
hypergraph.

Finally, we will compute out all the H-spectra of the class of
hyperstars, the hypercycle of size $3$ and the hyperpath of length
$3$.   This will be useful for research on the second smallest
Laplacian H-eigenvalue of $k$-uniform hypergraphs.

For discussion on the eigenvectors of the zero Laplacian and
signless Laplacian eigenvalues of a $k$-uniform hypergraph, see
\cite{hq13}.   For discussion on eigenvalues of adjacency tensors
and the other types of Laplacian tensors of $k$-uniform hypergraphs,
see \cite{cd12,hq12a,lqy12,pz12,r09,rp09,xc12c,xc12a} and references
therein.

The rest of this paper is organized as follows. Definitions on
eigenvalues of tensors and uniform hypergraphs are presented in the
next section. Cored hypergraphs and power hypergraphs are introduced
there. We discuss in Section \ref{sec-ch} some properties on the
cored hypergraphs. An even-uniform cored hypergraph has equality for
the largest Laplacian and the singless Laplacian H-eigenvalues.
Sunflowers are introduced and investigated in Section \ref{sec-sps}.
We compute out the largest Laplacian H-eigenvalues of even-uniform
sunflowers and prove that they are equal to the maximum degrees,
i.e., $2$, for odd-uniform sunflowers. We show in Section
\ref{sec-ph-1} that the largest Laplacian H-eigenvalues of
odd-uniform hypercycles and hyperpaths are equal to the maximum
degrees, i.e., $2$. We make a conjecture in Section \ref{sec-ph-2}
that the largest H-eigenvalues of even-uniform power hypergraphs
with respect to the same underlying usual graph are strictly
decreasing as $k$ increasing. This conjecture is proved to be true for hyperstars and
hypercycles. In Section \ref{sec-sph}, we compute out all the
H-eigenvalues of hyperstars, the hyperpath of length $3$ and the
hypercycle of size $3$. Some final remarks are made in the last
section.

\section{Preliminaries}\label{sec-p}
\setcounter{Theorem}{0} \setcounter{Proposition}{0}
\setcounter{Corollary}{0} \setcounter{Lemma}{0}
\setcounter{Definition}{0} \setcounter{Remark}{0}
\setcounter{Conjecture}{0}  \setcounter{Example}{0} \hspace{4mm}
\subsection{H-Eigenvalues of Tensors}\label{s-et}
In this subsection, some definitions of H-eigenvalues of tensors are presented. For comprehensive references, see \cite{q05,hhlq12} and references therein. Especially, for spectral hypergraph theory oriented facts on eigenvalues of tensors, please see \cite{q12a,hq13}.

Let $\mathbb R$ be the field of real numbers and $\mathbb R^n$ the
$n$-dimensional real space. $\mathbb R^n_+$ denotes the nonnegative
orthant of $\mathbb R^n$. For integers $k\geq 3$ and $n\geq 2$, a
real tensor $\mathcal T=(t_{i_1\ldots i_k})$ of order $k$ and
dimension $n$ refers to a multidimensional array (also called hypermatrix)
with entries $t_{i_1\ldots i_k}$ such that $t_{i_1\ldots
i_k}\in\mathbb{R}$ for all $i_j\in[n]:=\{1,\ldots,n\}$ and
$j\in[k]$. Tensors are always referred to $k$-th order real tensors
in this paper, and the dimensions will be clear from the content.
Given a vector $\mathbf{x}\in \mathbb{R}^{n}$, ${\cal
T}\mathbf{x}^{k-1}$ is defined as an $n$-dimensional vector such
that its $i$-th element being
$\sum\limits_{i_2,\ldots,i_k\in[n]}t_{ii_2\ldots i_k}x_{i_2}\cdots
x_{i_k}$ for all $i\in[n]$. Let ${\cal I}$ be the identity tensor of
appropriate dimension, e.g., $i_{i_1\ldots i_k}=1$ if and only if
$i_1=\cdots=i_k\in [n]$, and zero otherwise when the dimension is
$n$. The following definition was introduced by Qi \cite{q05}.
\begin{Definition}\label{def-00}
Let $\mathcal T$ be a $k$-th order $n$-dimensional real tensor. For
some $\lambda\in\mathbb{R}$, if polynomial system $\left(\lambda
{\cal I}-{\cal T}\right)\mathbf{x}^{k-1}=0$ has a solution
$\mathbf{x}\in\mathbb{R}^n\setminus\{0\}$, then $\lambda$ is called
an H-eigenvalue and $\mathbf x$ an
H-eigenvector.
\end{Definition}
H-eigenvalues are real numbers, by Definition \ref{def-00}. By
\cite{hhlq12,q05}, we have that the number of H-eigenvalues of a
real tensor is finite. By \cite{q12a}, we have that all the tensors
considered in this paper have at least one H-eigenvalue. Hence, we
can denote by $\lambda(\mathcal T)$ as the largest H-eigenvalue of a
real tensor $\mathcal T$.

For a subset $S\subseteq [n]$, we denoted by $|S|$ its cardinality, and $\mbox{sup}(\mathbf x):=\{i\in[n]\;|\;x_i\neq 0\}$ is the {\em support} of $\mathbf x$.

\subsection{Uniform Hypergraphs}

In this subsection, we present some essential notions of uniform hypergraphs which will be used in the sequel. Please refer to \cite{b73,c97,bh11,hq13,q12a} for comprehensive references.

In this paper, unless stated otherwise, a hypergraph means an undirected simple $k$-uniform hypergraph $G$ with vertex set $V$, which is labeled as $[n]=\{1,\ldots,n\}$, and edge set $E$.
By $k$-uniformity, we mean that for every edge $e\in E$, the cardinality $|e|$ of $e$ is equal to $k$. Throughout this paper, $k\geq 3$ and $n\geq k$. Moreover, since the trivial hypergraph (i.e., $E=\emptyset$) is of less interest, we consider only hypergraphs having at least one edge (i.e., nontrivial) in this paper.

For a subset $S\subset [n]$, we denote by $E_S$ the set of edges $\{e\in E\;|\;S\cap e\neq\emptyset\}$. For a vertex $i\in V$, we simplify $E_{\{i\}}$ as $E_i$. It is the set of edges containing the vertex $i$, i.e., $E_i:=\{e\in E\;|\;i\in e\}$. The cardinality $|E_i|$ of the set $E_i$ is defined as the {\em degree} of the vertex $i$, which is denoted by $d_i$. Two different vertices $i$ and $j$ are {\em connected} to each other (or the pair $i$ and $j$ is connected), if there is a sequence of edges $(e_1,\ldots,e_m)$ such that $i\in e_1$, $j\in e_m$ and $e_r\cap e_{r+1}\neq\emptyset$ for all $r\in[m-1]$. A hypergraph is called {\em connected}, if every pair of different vertices of $G$ is connected. Let $S\subseteq V$, the hypergraph with vertex set $S$ and edge set $\{e\in E\;|\;e\subseteq S\}$ is called the {\em sub-hypergraph} of $G$ induced by $S$. We will denote it by $G_S$. A hypergraph is {\em regular} if $d_1=\cdots=d_n=d$.
A hypergraph $G=(V,E)$ is {\em complete} if $E$ consists of all the possible edges. In this case, $G$ is regular, and moreover $d_1=\cdots=d_n=d={n-1\choose k-1}$.
In the sequel, unless stated otherwise, all the notations introduced above are reserved for the specific meanings.

 For the sake of simplicity, we mainly consider connected hypergraphs in the subsequent analysis. By the techniques in \cite{q12a,hq13}, the conclusions on connected hypergraphs can be easily generalized to general hypergraphs.

The following definition for the Laplacian tensor and signless Laplacian tensor was proposed by Qi \cite{q12a}.
\begin{Definition}\label{def-l}
Let $G=(V,E)$ be a $k$-uniform hypergraph. The {\em adjacency tensor} of $G$ is defined as the $k$-th order $n$-dimensional tensor $\mathcal A$ whose $(i_1 \ldots i_k)$-entry is:
\begin{eqnarray*}
a_{i_1 \ldots i_k}:=\left\{\begin{array}{cl}\frac{1}{(k-1)!}&if\;\{i_1,\ldots,i_k\}\in E,\\0&\mbox{otherwise}.\end{array}\right.
\end{eqnarray*}
Let $\mathcal D$ be a $k$-th order $n$-dimensional diagonal tensor with its diagonal element $d_{i\ldots i}$ being $d_i$, the degree of vertex $i$, for all $i\in [n]$. Then $\mathcal L:=\mathcal D-\mathcal A$ is the {\em Laplacian tensor} of the hypergraph $G$, and $\mathcal Q:=\mathcal D+\mathcal A$ is the {\em signless Laplacian tensor} of the hypergraph $G$.
\end{Definition}

By \cite{q12a}, zero is always the smallest H-eigenvalue of
$\mathcal L$ and $\mathcal Q$, and we have $\lambda(\mathcal L) \le
\lambda(\mathcal Q) \le 2d$, where $d$ is the maximum degree of $G$.

In the following, we introduce the class of cored hypergraphs.
\begin{Definition}\label{def-cor}
Let $G=(V,E)$ be a $k$-uniform hypergraph. If for every edge $e\in E$, there is a vertex $i_e\in e$ such that the degree of the vertex $i_e$ is one, then $G$ is a {\em cored hypergraph}. A vertex with degree being one is a {\em cored vertex}, and a vertex with degree being larger than one is an {\em intersectional vertex}.
\end{Definition}

Let $G=(V,E)$ be an ordinary graph. For every $k\geq 3$, we can
introduce a hypergraph by blowing up the edges of $G$.
\begin{Definition}\label{def-power}
Let $G=(V,E)$ be a $2$-uniform graph. For any $k\geq 3$, the $k$th
power of $G$, $G^k:=(V^k,E^k)$ is defined as the $k$-uniform
hypergraph with the set of edges being
$E^k:=\{e\cup\{i_{e,1},\ldots,i_{e,k-2}\}\;|\;e\in E\}$, and the set
of vertices being $V^k:=V\cup\{i_{e,1},\ldots,i_{e,k-2},\;e\in E\}$.
\end{Definition}
It is easy to see that the class of power hypergraphs is a subclass
of cored hypergraphs. The classes of hyperstars and hypercycles are
introduced in \cite{hqx13}. It can be seen that the classes of
hyperstars and hypercycles are subclasses of power hypergraphs.
Actually, a $k$-uniform hyperstar (respectively hypercycle) is the
$k$th power of a star (respectively cycle) graph.

We present in Figure 1 an example of an ordinary graph and its $3$rd
and $4$th power hypergraphs.
\begin{figure}[htbp]
\centering
\includegraphics[width=1.6in]{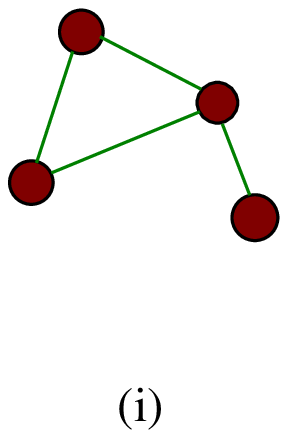}
\includegraphics[width=1.8in]{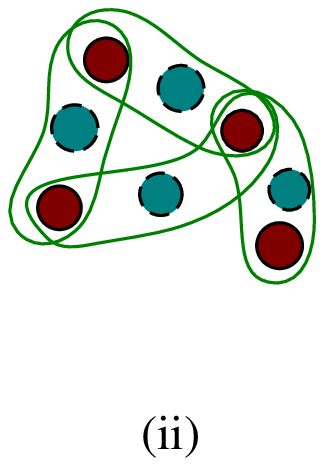}
\includegraphics[width=2.0in]{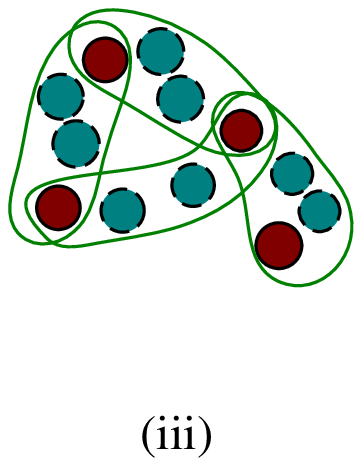}
\caption{(i) is an example of an usual graph. (ii) and (iii) are respectively the $3$-uniform and $4$-uniform power hypergraphs of the graph in (i). A solid disk represents a vertex. For the usual graph, a line connecting two points is an edge. For hypergraphs, an edge is pictured as a closed curve with the containing solid disks the vertices in that edge. The newly added vertices are in different color (also in dashed margins). }
\end{figure}

For completeness, we include the defintiions for hyperstars and hypercycles in Definitions \ref{def-bi-hm} and \ref{def-hc} respectively.

\begin{Definition}\label{def-bi-hm}
Let $G=(V,E)$ be a $k$-uniform hypergraph. If there is a disjoint partition of the vertex set $V$ as $V=V_0\cup V_1\cup\cdots\cup V_d$ such that $|V_0|=1$ and $|V_1|=\cdots=|V_d|=k-1$, and $E=\{V_0\cup V_i\;|\;i\in [d]\}$, then $G$ is called a {\em hyperstar}. The degree $d$ of the vertex in $V_0$, which is called the {\em heart}, is the {\em size} of the hyperstar. The edges of $G$ are {\em leaves}, and the vertices other than the heart are vertices of leaves.
\end{Definition}

\begin{Definition}\label{def-hc}
Let $G=(V,E)$ be a $k$-uniform nontrivial
hypergraph. If we can number the vertex set $V$ as $V:=\{i_{1,1},\ldots,i_{1,k-1},\ldots,i_{d,1},\ldots,i_{d,k-1}\}$ for some positive integer $d$ such that $E=\{\{i_{1,1},\ldots,i_{1,k-1},i_{2,1}\},\{i_{2,1},\ldots,i_{2,k-1},i_{3,1}\},\ldots,\{i_{d-1,1},\ldots,i_{d-1,k-1},i_{d,1}\},\{i_{d,1},\ldots,i_{d,k-1},i_{1,1}\}\}$,
then $G$ is called a {\em hypercycle}. $d$ is the {\em size} of the hypercycle.
\end{Definition}
It is easy to see that a $k$-uniform hyperstar of size $s>0$ has $n=s(k-1)+1$ vertices, a $k$-uniform hypercycle of size $s>0$ has $n=s(k-1)$ vertices, and they are both connected.

Besides hyperstars and hypercycles, power hypergraphs contain hyperpaths. Hyperpaths are power hypergraphs of usual paths. We state it in the next definition.

\begin{Definition}\label{def-path}
Let $G=(V,E)$ be a $k$-uniform hypergraph. If we can number the vertex set $V$ as $V:=\{i_{1,1},\ldots,i_{1,k},i_{2,2},\ldots,i_{2,k},\ldots,i_{d-1,2},\ldots,i_{d-1,k},i_{d,2},\ldots,i_{d,k}\}$ for some positive integer $d$ such that $E=\{\{i_{1,1},\ldots,i_{1,k}\},\{i_{1,k},i_{2,2},\ldots,i_{2,k}\},\ldots,\{i_{d-1,k},i_{d,2},\ldots,i_{d,k}\}\}$, then $G$ is a {\em hyperpath}. $d$ is the {\em length} of the hyperpath.
\end{Definition}
Figure 2 is an axample of a $3$-uniform hyperpath.
\begin{figure}[htbp]
\centering
\includegraphics[width=2.6in]{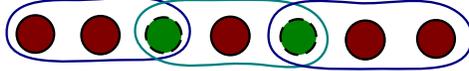}
\caption{An example of a $3$-uniform hyperpath of length $3$. The intersectional vertices are in red (also in dashed margins).}
\end{figure}

The notions of odd-bipartite and even-bipartite even-uniform hypergraphs are introduced in \cite{hq13b}.
\begin{Definition}\label{def-bi-odd}
Let $k$ be even and $G=(V,E)$ be a $k$-uniform
hypergraph. It is called {\em odd-bipartite} if either it is trivial
(i.e., $E=\emptyset$) or there is a disjoint partition of the vertex
set $V$ as $V=V_1\cup V_2$ such that $V_1,V_2\neq \emptyset$ and
every edge in $E$ intersects $V_1$ with exactly an odd number of
vertices.
\end{Definition}

\section{Cored Hypergraphs}\label{sec-ch}
\setcounter{Theorem}{0} \setcounter{Proposition}{0}
\setcounter{Corollary}{0} \setcounter{Lemma}{0}
\setcounter{Definition}{0} \setcounter{Remark}{0}
\setcounter{Conjecture}{0}  \setcounter{Example}{0} \hspace{4mm}

Some facts on the H-eigenvalues and H-eigenvectors of the Laplacian tensor of a cored hypergraph is discussed in this section.

\subsection{General Cases}
In this subsection, we establish some facts that all cored hypergraphs share.

The next lemma says that their H-eigenvectors have special
structures.
\begin{Lemma}\label{lem-ch-1}
Let $G=(V,E)$ be a $k$-uniform cored hypergraph and $\mathbf x\in\mathbb R^n$ be an H-eigenvalue of its Laplacian tensor $\mathcal L$ corresponding to an H-eigenvalue $\lambda\neq 1$. If there are two cored vertices $i,j$ in an edge $e\in E$, then $|x_i|=|x_j|$. Moreover, $x_i=x_j$ when $k$ is an odd number.
\end{Lemma}

\noindent {\bf Proof.} By the definition of H-eigenvalues and the fact that $i$ and $j$ are cored vertices, we have
\begin{eqnarray*}
\lambda x_i^{k-1}=(\mathcal L\mathbf x^{k-1})_i=x_i^{k-1}-x_j\prod_{s\in e\setminus\{i,j\}}x_s,\;\;\lambda x_j^{k-1}=(\mathcal L\mathbf x^{k-1})_j=x_j^{k-1}-x_i\prod_{s\in e\setminus\{i,j\}}x_s.
\end{eqnarray*}
Hence,
\begin{eqnarray*}
(\lambda-1) x_i^{k}=(\lambda-1) x_j^k.
\end{eqnarray*}
Since $\lambda\neq 1$, we have that $|x_i|=|x_j|$. Moreover, when $k$ is odd, we see that $x_i=x_j$.  \ep

By \cite[Theorem 4]{q12a}, we have the following lemma.
\begin{Lemma}\label{lem-ch-2}
Let $G=(V,E)$ be a $k$-uniform hypergraph with its maximum degree
$d>0$ and $\mathcal L$ be its Laplacian
tensor. Then $\lambda (\mathcal L)\geq d$.
\end{Lemma}

\begin{Lemma}\label{lem-ch-3}
Let $G=(V,E)$ be a $k$-uniform cored hypergraph and $\mathbf x\in\mathbb R^n$ be an H-eigenvalue of its Laplacian tensor $\mathcal L$ corresponding to $\lambda \geq 1$. Then, $\prod_{s\in e}x_s\leq 0$ for all $e\in E$ when $k$ is even; and $\prod_{s\in e\setminus\{i_e\}}x_s\leq 0$ for all $e\in E$ when $k$ is odd. Here $i_e\in e$ is a cored vertex.
\end{Lemma}

\noindent {\bf Proof.} Suppose that $i$ is a cored vertex of an arbitrary but fixed edge $e\in E$. If $\lambda=1$, then
\begin{eqnarray*}
 x_i^{k-1}=\lambda x_i^{k-1}=x_i^{k-1}-\prod_{s\in e\setminus\{i\}}x_s
\end{eqnarray*}
implies that $\prod_{s\in e\setminus\{i\}}x_s=0$. We are done.

In the following, suppose that $\lambda>1$. Then,
\begin{eqnarray*}
(\lambda -1)x_i^{k-1}=-\prod_{s\in e\setminus\{i\}}x_s.
\end{eqnarray*}
When $k$ is odd, $x_i^{k-1}\geq 0$. Then the result follows.
When $k$ is even, we have $\prod_{s\in e}x_s\leq 0$ since $(\lambda -1)x_i^k=-\prod_{s\in e}x_s$. \ep

By Lemmas \ref{lem-ch-2} and \ref{lem-ch-3}, we get the next proposition.
\begin{Proposition}\label{prop-ch-0}
Let $G=(V,E)$ be a $k$-uniform cored hypergraph and $\mathbf x\in\mathbb R^n$ be an H-eigenvalue of its Laplacian tensor $\mathcal L$ corresponding to $\lambda(\mathcal L)$. Then, $\prod_{s\in e}x_s\leq 0$ for all $e\in E$ when $k$ is even; and $\prod_{s\in e\setminus\{i_e\}}x_s\leq 0$ for all $e\in E$ when $k$ is odd. Here $i_e\in e$ is a cored vertex.
\end{Proposition}

By \cite[Theorem 5.1]{hqx13}, we can get the next proposition.
\begin{Proposition}\label{prop-ch-1}
Let $k$ be even and $G=(V,E)$ be a $k$-uniform cored hypergraph. Let $\mathcal L$ and $\mathcal Q$ be the Laplacian tensor and signless Laplacian tensor of $G$ respectively. Then $G$ is odd-bipartite, and hence $\lambda(\mathcal L)=\lambda(\mathcal Q)$.
\end{Proposition}

\noindent {\bf Proof.} For all $e\in E$, let $i_e\in e$ be a cored vertex. Set $V_1:=\{i_e\;|\;e\in E\}$ and $V_2:=V\setminus V_1$. Then it is easy to see that $V=V_1\cup V_2$ is an odd-bipartition (Definition \ref{def-bi-odd}). Thus, the result follows from \cite[Theorem 5.1]{hqx13}.  \ep

Actually, we can get the next proposition.
\begin{Proposition}\label{prop-ch-2}
Let $k$ be even and $G=(V,E)$ be a $k$-uniform cored hypergraph. Let $\mathcal L$ and $\mathcal Q$ be the Laplacian tensor and signless Laplacian tensor of $G$ respectively. For every $e\in E$, let $i_e\in e$ be a cored vertex.
\begin{itemize}
\item [(i)] If $\mathbf x\in\mathbb R_+^n$ is an H-eigenvector of $\mathcal Q$ corresponding to $\lambda(\mathcal Q)$, then $\mathbf y\in\mathbb R^n$ is an H-eigenvector of $\mathcal L$ corresponding to $\lambda(\mathcal L)$ with $y_{i_e}=-x_{i_e}$ for all $e\in E$ and $y_j=x_j$ for the others.
\item [(ii)] If $\mathbf x\in\mathbb R^n$ is an H-eigenvector of $\mathcal L$ corresponding to $\lambda(\mathcal L)$, then $\mathbf y\in\mathbb R_+^n$ is an H-eigenvector of $\mathcal Q$ corresponding to $\lambda(\mathcal Q)$ with $y_i=|x_i|$ for all $j\in [n]$.
\end{itemize}
\end{Proposition}

\noindent {\bf Proof.} The results follow from Definition \ref{def-00} and Proposition \ref{prop-ch-0}. \ep

\subsection{Sunflowers}\label{sec-sps}

Obviously, not all cored hypergraphs are power hypergraphs. Among the others, the class of sunflowers is investigated.
\begin{Definition}\label{def-sf}
Let $G=(V,E)$ be a $k$-uniform hypergraph. If we can number the vertex set $V$ as $V:=\{i_{1,1},\ldots,i_{1,k},\ldots,i_{k-1,1},\ldots,i_{k-1,k},i_k\}$ such that the set of edges being $E=\{\{i_{1,1},\ldots,i_{1,k}\},\ldots,\{i_{k-1,1},\ldots,i_{k-1,k}\},\{i_{1,1},\ldots,i_{k-1,1},i_k\}\}$, then $G$ is a {\em sunflower}.
\end{Definition}
Note that the sunflower for every positive integer $k$ is unique, in the sense that by a possible renumbering the vertices two $k$-uniform sunflowers are the same.
Figure 3 is an example of the $4$-uniform sunflower.
\begin{figure}[htbp]
\centering
\includegraphics[width=2.6in]{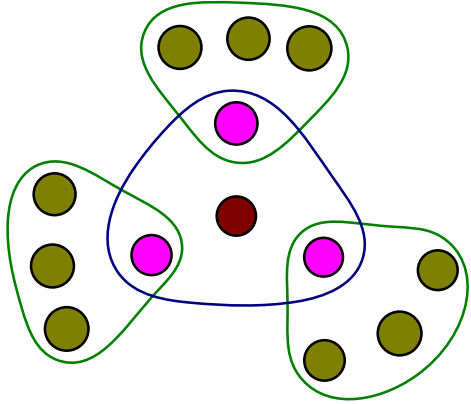}
\caption{ An example of the $4$-uniform sunflower. }
\end{figure}

The next proposition finds out the largest H-eigenvalue of the Laplacian tensor of an even-uniform sunflower.
\begin{Proposition}\label{prop-ch-4}
Let $k$ be even and $G=(V,E)$ be the $k$-uniform sunflower. Let $\mathcal L$ be the Laplacian tensor of $G$. Then $G$ is odd-bipartite and $\lambda(\mathcal L)$ is the unique root of $(\mu-2)-\left(\frac{1}{\mu-1}\right)^{\frac{1}{k-1}}-\left(\frac{1}{\mu-1}\right)^{k-1}=0$ in the interval $(2,4)$.
\end{Proposition}

\noindent {\bf Proof.} Suppose that $V=\{i_{1,1},\ldots,i_{1,k},\ldots,i_{k-1,1},\ldots,i_{k-1,k},i_k\}$, and the set of edges is $E=\{\{i_{1,1},\ldots,i_{1,k}\},\ldots,\{i_{k-1,1},\ldots,i_{k-1,k}\},\{i_{1,1},\ldots,i_{k-1,1},i_k\}\}$.
Let $\mathcal Q$ be the signless Laplacian tensor of $G$. By Proposition \ref{prop-ch-1}, $G$ is odd-bipartite and $\lambda(\mathcal L)=\lambda(\mathcal Q)$.

By \cite[Theorem 3.20]{yy11}, \cite[Theorem 4]{q12b} and \cite[Lemma 3.1]{pz12} (see also \cite[Lemmas 2.2 and 2.3]{hq13}), if we can find a positive H-eigenvector $\mathbf x\in\mathbb R^n$ of $\mathcal Q$ corresponding to an H-eigenvalue $\mu$, then $\mu=\lambda(\mathcal Q)$.

Let $x_{i_k}=\alpha>0$, $x_{i_{j,1}}=1$, and
$x_{i_{j,2}}=\cdots=x_{i_{j,k}}=\gamma>0$ for all $j\in [k-1]$. Suppose that $\mathbf x$ is an H-eigenvector of $\mathcal Q$ corresponding to the H-eigenvalue $\mu=\lambda(\mathcal Q)$.
By Definition \ref{def-00}, we have
\begin{eqnarray*}
(\mu-1)\alpha^{k-1}=1,\;(\mu-2)=\alpha+\gamma^{k-1},\;\mbox{and}\;(\mu-1)\gamma^{k-1}=\gamma^{k-2}.
\end{eqnarray*}
By Lemma \ref{lem-ch-2}, we have $\mu\geq 2$. Thus, the first and the third equalities imply that $\alpha^{k-1}=\gamma$.
Hence,
\begin{eqnarray*}
(\mu-2)=\left(\frac{1}{\mu-1}\right)^{\frac{1}{k-1}}+\left(\frac{1}{\mu-1}\right)^{k-1}.
\end{eqnarray*}
Let $f(\mu):=(\mu-2)-\left(\frac{1}{\mu-1}\right)^{\frac{1}{k-1}}-\left(\frac{1}{\mu-1}\right)^{k-1}$. We have that $f(2)=-2<0$ and
\begin{eqnarray*}
f(4)=2-\left(\frac{1}{3}\right)^{\frac{1}{k-1}}-\frac{1}{3^{k-1}}>0.
\end{eqnarray*}
Thus, $f(\mu)=0$ does have a root in the interval $(2,4)$. Since $\mathcal Q$ has a unique positive H-eigenvector (\cite[Lemmas 2.2 and 2.3]{hq13}), the equation $(\mu-2)-\left(\frac{1}{\mu-1}\right)^{\frac{1}{k-1}}-\left(\frac{1}{\mu-1}\right)^{k-1}=0$ has a unique positive solution which is in the interval $(2,4)$. Hence, the result follows. \ep

The next proposition says that the largest H-eigenvalue of the Laplacian tensor of an odd-uniform sunflower is equal to the maximum degree, ie., $2$.
\begin{Proposition}\label{prop-ch-5}
Let $k$ be odd and $G=(V,E)$ be the $k$-uniform sunflower. Let $\mathcal L$ be the Laplacian tensor of $G$. Then $\lambda(\mathcal L)=2$.
\end{Proposition}

\noindent {\bf Proof.} Suppose that $V=\{i_{1,1},\ldots,i_{1,k},\ldots,i_{k-1,1},\ldots,i_{k-1,k},i_k\}$, and the set of edges is $E=\{\{i_{1,1},\ldots,i_{1,k}\},\ldots,\{i_{k-1,1},\ldots,i_{k-1,k}\},\{i_{1,1},\ldots,i_{k-1,1},i_k\}\}$. Let $\mathbf w\in\mathbb R^n$ be an H-eigenvector of $\mathcal L$ corresponding to $\lambda(\mathcal L)$. Then, we have
\begin{eqnarray*}
(\lambda(\mathcal L)-1)w_{i_{j,s}}^{k-1}=-w_{i_{j,1}}\prod_{t\in\{2,\ldots,k\}\setminus\{s\}}w_{i_{j,t}},\;\forall j\in[k-1],\;s\in\{2,\ldots,k\}.
\end{eqnarray*}
Thus, $(\lambda(\mathcal L)-1)w_{i_{j,s}}^k=-w_{i_{j,1}}\prod_{t\in\{2,\ldots,k\}}w_{i_{j,t}}$. Hence, $w_{i_{i,2}}=\cdots=w_{i_{j,k}}=:z_j$ for all $j\in [k-1]$, since $\lambda(\mathcal L)\geq 2$ by Lemma \ref{lem-ch-2}. Let $x:=w_{i_k}$, and $y_j:=w_{i_{j,1}}$ for all $j\in [k-1]$.  Then, we have $y_j=(1-\lambda(\mathcal L))z_j$ if $z_j\neq 0$, for all $j\in [k-1]$. Moreover, by Definition \ref{def-00}, we have
\begin{eqnarray*}
(\lambda(\mathcal L)-2)y_j^{k-1}=-x\prod_{s\in[k-1]\setminus\{j\}}y_s-z_j^{k-1},\;\forall j\in [k-1].
\end{eqnarray*}
Consequently, we have
\begin{eqnarray*}
\left[(\lambda(\mathcal L)-2)(1-\lambda(\mathcal L))^{k-1}+(1-\lambda(\mathcal L))\right] z_j^k=-x\prod_{s\in[k-1]}y_s,\;\forall j\in [k-1].
\end{eqnarray*}
Thus, either $z_1=\cdots=z_{k-1}=:z$ since $k$ is odd, or $(\lambda(\mathcal L)-2)(1-\lambda(\mathcal L))^{k-1}+(1-\lambda(\mathcal L))=0$ and $x\prod_{s\in[k-1]}y_s=0$.

If $(\lambda(\mathcal L)-2)(1-\lambda(\mathcal L))^{k-1}+(1-\lambda(\mathcal L))=0$, then $(\lambda(\mathcal L)-2)(1-\lambda(\mathcal L))^{k-2}+1=0$ since $\lambda(\mathcal L)\geq 2$. We must have that $\lambda(\mathcal L)>2$, since $f(2)=1>0$ and $f$ is decreasing in $(2,\infty)$. Here $f(\mu):=(\mu-2)(1-\mu)^{k-2}+1$. On the other hand, if all $y_t$ with $t\in [k-1]$ are zero, we get that $\mathbf w=0$ which is a contradiction. Hence, we must have that $x\prod_{s\in[k-1]\setminus\{j\}}y_s=0$ and $y_j\neq 0$ for some $j\in [k-1]$, since $x\prod_{s\in[k-1]}y_s=0$. These two facts will contradict the fact that
\begin{eqnarray*}
(\lambda(\mathcal L)-2)y_j^{k-1}=-x\prod_{s\in[k-1]\setminus\{j\}}y_s-z_j^{k-1}=-z_j^{k-1},
\end{eqnarray*}
since $k-1$ is even. Thus, this situation can never happen.

If $z_1=\cdots=z_{k-1}=z\neq 0$, then $y_1=\cdots=y_{k-1}=:y\neq 0$ since $y_j=(1-\lambda(\mathcal L))z_j$.
By Definition \ref{def-00}, we have $0\leq(\lambda(\mathcal L)-1)x^{k-1}=-y^{k-1}\leq 0$. Hence $x=y=0$, which is a contradiction. Then, we must have that $z=0$. The equations of the H-eigenvalue $\lambda(\mathcal L)$ become
\begin{eqnarray*}
(\lambda(\mathcal L)-1)x^{k-1}=-\prod_{s\in[k-1]}y_s,\;(\lambda(\mathcal L)-2)y_j^{k-1}=-x\prod_{s\in[k-1]\setminus\{j\}}y_s,\;\forall j\in [k-1].
\end{eqnarray*}
If $\lambda(\mathcal L)>2$, we must have $y_1=\cdots=y_{k-1}=:y$, since $(\lambda(\mathcal L)-2)y_j^k=-x\prod_{s\in[k-1]}y_s$ for all $j\in [k-1]$. Similarly, we must have $x=y=0$ which is a contradiction. Hence, $\lambda(\mathcal L)=2$. An H-eigenvector would be $y_1=1$ and the rest are zero.
\ep

\section{Power Hypergraphs}\label{sec-ph}
\setcounter{Theorem}{0} \setcounter{Proposition}{0}
\setcounter{Corollary}{0} \setcounter{Lemma}{0}
\setcounter{Definition}{0} \setcounter{Remark}{0}
\setcounter{Conjecture}{0}  \setcounter{Example}{0} \hspace{4mm}
Some facts on the H-eigenvalues and H-eigenvectors of the Laplacian tensor of a power hypergraph are investigated in this section.

\subsection{Odd-Uniform Power Hypergraphs}\label{sec-ph-1}

In this subsection, we show that the largest H-eigenvalue of the Laplacian tensor of an odd-uniform hypercycle (hyperpath) is equal to the maximum degree, i.e., $2$.

The next lemma is useful.
\begin{Lemma}\label{lem-ph-1}
Let $k$ be odd, $G=(V,E)$ be a $k$-uniform power hypergraph and $\mathbf x\in\mathbb R^n$ be an H-eigenvalue of its Laplacian tensor $\mathcal L$ corresponding to $\lambda\neq 1$. Let $e\in E$ be an arbitrary but fixed edge.
\begin{itemize}
\item [(i)] If $e$ has only one intersectional vertex $i$, and $x_s\neq 0$ for some cored vertex $s\in e$, then $(1-\lambda)x_s=x_i$.
\item [(ii)] If $e$ has two intersectional vertices $i$ and $j$, and $x_s\neq 0$ for some cored vertex $s\in e$, then $x_ix_j=(1-\lambda)x_s^2$.
\end{itemize}
\end{Lemma}

\noindent {\bf Proof.} For (i), by Definition \ref{def-00} and Lemma \ref{lem-ch-1}, we have
\begin{eqnarray*}
\lambda x_s^{k-1}=x_s^{k-1}-x_s^{k-2}x_i.
\end{eqnarray*}
Thus, $x_i=(1-\lambda)x_s$.

For (ii), by Definition \ref{def-00} and Lemma \ref{lem-ch-1}, we have
\begin{eqnarray*}
\lambda x_s^{k-1}=x_s^{k-1}-x_s^{k-3}x_ix_j.
\end{eqnarray*}
Thus, $x_ix_j=(1-\lambda)x_s^2$.
\ep

The next corollary is a direct consequence of Lemma \ref{lem-ph-1}.
\begin{Corollary}\label{cor-ph-1}
Let $k$ be odd and $G=(V,E)$ be a $k$-uniform power hypergraph and $\mathbf x\in\mathbb R^n$ be an H-eigenvalue of its Laplacian tensor $\mathcal L$ corresponding to $\lambda>1$. Let $e\in E$ be an arbitrary but fixed edge.
\begin{itemize}
\item [(i)] If $e$ has only one intersectional vertex $i$, and $x_s\neq 0$ for some cored vertex $s\in e$, then $x_ix_s<0$.
\item [(ii)] If $e$ has two intersectional vertices $i$ and $j$, and $x_s\neq 0$ for some cored vertex $s\in e$, then $x_ix_j<0$.
\end{itemize}
\end{Corollary}

The next proposition is for hypercycle.
\begin{Proposition}\label{prop-ph-2}
Let $k$ be odd and $G=(V,E)$ be a $k$-uniform hypercycle with size being $r\geq 2$. Let $\mathcal L$ be its Laplacian tensor. Then $\lambda(\mathcal L)=2$.
\end{Proposition}

\noindent {\bf Proof.} Suppose that
\begin{eqnarray*}
E&=&\{\{j_{1},i_{1,1}\ldots,i_{1,k-2},j_2\},\{j_2,i_{2,1}\ldots,i_{2,k-2},j_3\},\ldots,\\
&&\{j_{r-1},i_{r-1,1}\ldots,i_{r-1,k-2},j_r\},\{j_r,i_{r,1}\ldots,i_{r,k-2},j_1\}\}.
\end{eqnarray*}
Let $\mathbf z\in\mathbb R^{r(k-1)}$ be an H-eigenvector of $\mathcal L$ corresponding to $\lambda(\mathcal L)$. Let $y_s:=z_{j_s}$ for $s\in [r]$. By Lemma \ref{lem-ch-1}, we have that $z_{i_{s,1}}=\cdots=z_{i_{s,k-2}}$. Let $x_s:=z_{i_{s,1}}$ for $s\in [r]$.

(I). If $x_s\neq 0$ for all $s\in [r]$, by Lemma \ref{lem-ph-1}, we have that $y_sy_{s+1}<0$ with $y_{r+1}:=y_1$ for all $s\in [r]$.
Thus, if $r$ is odd, we get a contradiction by the rule of signs. In the following, we assume that $r$ is even and $\lambda(\mathcal L)>2$.
By Definition \ref{def-00}, we have
\begin{eqnarray*}
(\lambda(\mathcal L)-2)y_s^k=-x_{s-1}^{k-2}y_{s-1}y_s-x_s^{k-2}y_sy_{s+1},\;\forall s\in [r]
\end{eqnarray*}
with the convention that $x_0=x_r$, $x_{r+1}=x_1$, and $y_0=y_r$, $y_{r+1}=y_1$. Thus, we have
\begin{eqnarray*}
y_1^k-y_2^k+y_3^k+\cdots+y_{r-1}^k-y_r^k=0.
\end{eqnarray*}
Since $y_sy_{s+1}<0$ and $k$ is odd, we get a contradiction, since $y_1^k-y_2^k+y_3^k+\cdots+y_{r-1}^k-y_r^k$ should be positive (respectively negative) if we let $y_1$ be positive (respectively negative). Consequently, $\lambda(\mathcal L)=2$.

(II). Suppose that $x_s=0$ for some $s\in [r]$. If all $x_s=0$ for $s\in [r]$. Then, we are done by Definition \ref{def-00}. In the following, we assume that $\lambda(\mathcal L)>2$. By Lemma \ref{lem-ph-1}, we have that $y_sy_{s+1}<0$ whenever $x_s\neq 0$.

Without loss of generality, we assume that $x_r=0$, $y_1\neq 0$, $x_1\neq 0$, $\cdots$, $y_m\neq 0$ for some $m\leq r$. Moreover, we can assume that $y_1>0$. By Definition \ref{def-00}, we have
\begin{eqnarray*}
(\lambda(\mathcal L)-2)y_1^{k-1}=-x_1^{k-2}y_2.
\end{eqnarray*}
Since $y_2<0$, we have that $x_1>0$. Also,
\begin{eqnarray*}
(\lambda(\mathcal L)-2)y_2^{k-1}=-x_1^{k-2}y_1-x_2^{k-2}y_3.
\end{eqnarray*}
Then, we must have $x_2<0$. Inductively, we have $x_sy_{s+1}<0$. Hence, $x_{m-1}y_m<0$ which implies that $x_{m-1}y_{m-1}>0$. By Definition \ref{def-00}, we have
\begin{eqnarray*}
(\lambda(\mathcal L)-2)y_m^{k-1}=-x_{m-1}^{k-2}y_{m-1}<0.
\end{eqnarray*}
Thus, a contradiction is derived. Consequently, $\lambda(\mathcal L)=2$. \ep

The next proposition is on hyperpath.
\begin{Proposition}\label{prop-ph-3}
Let $k$ be odd and $G=(V,E)$ be a $k$-uniform hyperpath with length being $r\geq 3$. Let $\mathcal L$ be its Laplacian tensor. Then $\lambda(\mathcal L)=2$.
\end{Proposition}

\noindent {\bf Proof.} Suppose that
\begin{eqnarray*}
E&=&\{\{i_{1,1}\ldots,i_{1,k-1},j_1\},\{j_1,i_{2,1}\ldots,i_{2,k-2},j_2\},\ldots,\\
&&\{j_{r-2},i_{r-1,1}\ldots,i_{r-1,k-2},j_{r-1}\},\{j_{r-1},i_{r,1}\ldots,i_{r,k-1}\}\}.
\end{eqnarray*}
Let $\mathbf z\in\mathbb R^{r(k-1)+1}$ be an H-eigenvector of $\mathcal L$ corresponding to $\lambda(\mathcal L)$. Let $y_s:=z_{j_s}$ for $s\in [r-1]$. By Lemma \ref{lem-ch-1}, we have that $z_{i_{s,1}}=\cdots=z_{i_{s,k-2}}$ for $s\in\{2,\ldots,r-1\}$, and $z_{i_{i,1}}=\cdots=z_{i_{1,k-1}}$, $z_{i_{r,1}}=\cdots=z_{i_{r,k-1}}$. Let $x_s:=z_{i{s,1}}$ for $s\in [r]$. We assume that $\lambda(\mathcal L)>2$ and derive a contradiction case by case. The proof is in the same spirit of and similar to that for Proposition \ref{prop-ph-2}, we include it here for completeness.

(I). If both $x_1$ and $x_r$ are zero, then the proof is same as that for the proof (II) in Proposition \ref{prop-ph-2}, since we always has a piece of the hyperpath with $y_{t},x_{t+1},\ldots,y_m\neq 0$ for some $m\geq t\geq 1$.

(II). If $x_1\neq 0$ and $x_r=0$, then we can find some $m\geq 1$ such that $x_1,y_1,\ldots,y_m\neq 0$. We can assume that $y_1>0$. By Lemma \ref{lem-ph-1}, we have $x_1<0$.
By Definition \ref{def-00}, we have
\begin{eqnarray*}
(\lambda(\mathcal L)-2)y_1^{k-1}=-x_1^{k-1}-x_2^{k-2}y_2.
\end{eqnarray*}
Since $y_2<0$ by Lemma \ref{lem-ph-1}, we have that $x_2>0$. Also,
\begin{eqnarray*}
(\lambda(\mathcal L)-2)y_2^{k-1}=-x_2^{k-2}y_1-x_3^{k-2}y_3.
\end{eqnarray*}
Then, we must have $x_3>0$. Inductively, we have $x_sy_{s}<0$. Hence, $x_{m}y_m<0$ which implies that $x_{m}y_{m-1}>0$. By Definition \ref{def-00}, we have
\begin{eqnarray*}
(\lambda(\mathcal L)-2)y_m^{k-1}=-x_{m}^{k-2}y_{m-1}<0.
\end{eqnarray*}
Thus, a contradiction is derived. Consequently, $\lambda(\mathcal L)=2$.

(III). The proof for the case $x_1=0$ and $x_r\neq 0$ is similar. Actually, it follows from (II) immediately by renumbering the indices.

(IV). If $x_s\neq 0$ for all $s\in [r]$. Similar to (II), we have that $x_sy_s<0$. Particularly, we have $x_{r-1}y_{r-1}<0$ which implies that $x_{r-1}y_{r-2}>0$. By Definition \ref{def-00}, we have
\begin{eqnarray*}
(\lambda(\mathcal L)-2)y_{r-1}^{k-1}=-x_{r-1}^{k-2}y_{r-2}-x_r^{k-1}<0.
\end{eqnarray*}
Thus, a contradiction is derived. Consequently, $\lambda(\mathcal L)=2$.

The other cases can be handled by the above proof as well. \ep

By Propositions \ref{prop-sph-2}, \ref{prop-ph-2} and \ref{prop-ph-3}, we get that \cite[Conjecture 3.1]{hqx13} has a negative answer.

\subsection{Even-Uniform Power Hypergraphs}\label{sec-ph-2}

We have a conjecture for even-uniform hypergraphs.
\begin{Conjecture}\label{thm-even-1}
Let $G=(V,E)$ be an usual graph, $k=2r$ be even and $G^k=(V^k,E^k)$ be the $k$-power hypergraph of $G$. Let $\mathcal L^k$ and $\mathcal Q^k$ be the Laplacian and signless Laplacian tensors of $G^k$ respectively. Then $\{\lambda(\mathcal L^k)=\lambda(\mathcal Q^k)\}$ is a strictly decreasing sequence.
\end{Conjecture}

By \cite[Theorem 3.1 and Corollary 5.1]{hqx13}, we have the next
proposition.

\begin{Proposition}\label{prop-ph-1}
Conjecture \ref{thm-even-1} is true for hyperstars and hypercycles.
\end{Proposition}

\noindent {\bf Proof.} For the case of hyperstars, by \cite[Theorem 3.1]{hqx13}, we have that $\lambda(\mathcal L^k)$ is the unique root of
\begin{eqnarray*}
(1-\mu)^{k-1}(\lambda-d)+d=0.
\end{eqnarray*}
Here $d$ is the size of the hyperstar. Let $f_k(\mu):=(1-\mu)^{k-1}(\lambda-d)+d$. We have $f_(k+1)(d)=d>0$ and
\begin{eqnarray*}
f_{k+1}(\lambda(\mathcal L^k))=(1-\lambda(\mathcal L^k))^k(\lambda(\mathcal L^k)-d)+d<(1-\lambda(\mathcal L^k))^{k-1}(\lambda(\mathcal L^k)-d)+d=0.
\end{eqnarray*}
Hence, $\lambda(\mathcal L^{k+1})\in (d, \lambda(\mathcal L^k))$.

The proof for the other case is similar. \ep

\section{H-Spectra of Special Power Hypergraphs}\label{sec-sph}
\setcounter{Theorem}{0} \setcounter{Proposition}{0}
\setcounter{Corollary}{0} \setcounter{Lemma}{0}
\setcounter{Definition}{0} \setcounter{Remark}{0}
\setcounter{Conjecture}{0}  \setcounter{Example}{0} \hspace{4mm}
We compute out all the H-eigenvalues of some special power hypergraphs in this section.

\subsection{Hyperstars}\label{sec-hs}
Let $G=(V,E)$ be a $k$-uniform hyperstar with $k\ge 3$ and the size $d\ge 2$, where $V=[n]$, $E=\{e_1,\cdots,e_d\}$, and $d_1=d$ (i.e., the vertex 1 is the heart). Let $\mathcal {L}= \mathcal {D}-\mathcal {A}$ be the Laplacian tensor of $G$. Then it is easy to see that the eigenvalue equations $(\lambda \mathcal {I}-\mathcal {L})\mathbf x^{k-1}=0$ are equivalent to the following set of relations:

\begin{eqnarray}\label{hyst-1}
(\lambda -d)x_1^{k-1}= - \sum_{i=1}^d \prod_{s\in e_i\setminus\{1\}}x_s,
\end{eqnarray}
and
\begin{eqnarray}\label{hyst-2}
(\lambda -1)x_j^{k-1}= -  \prod_{s\in e(j)\setminus \{j\}}x_s\;\; \forall j\in\{2,\ldots,n\},
\end{eqnarray}
where $e(j)$ denotes the unique edge containing the vertex $j$ for $j\ge 2$.

The next lemma strengthens Lemma \ref{lem-ph-1} for the case of hyperstars.

\begin{Lemma}\label{lem-0}
Let $G,k,d$ and $\mathcal {L}$ be as above. Suppose that $(\lambda, x)$ is an H-eigenpair of $\mathcal {L}$ with $\lambda \ne 1$. Then we have:
\begin{itemize}
\item [(i)] If $i,j\ge 2$ and $i,j$ are adjacent (there is an edge containing both $i$ and $j$), then $x_i=x_j$ when $k$ is odd, and  $|x_i|=|x_j|$ when $k$ is even.
\item [(ii)] If $i,j\ge 2$ and $x_i$, $x_j$ are both nonzero, then $x_i=x_j$ when $k$ is odd, and  $|x_i|=|x_j|$ when $k$ is even.
\end{itemize}
\end{Lemma}

\noindent {\bf Proof}: (i) follows from Lemma \ref{lem-ph-1}.

(ii) {\bf Case 1: $k$ is odd}. If $j\ge 2$ and $x_j\ne 0$, by \reff{hyst-2} and the result (i) of this lemma we also have $x_1= (1-\lambda )x_j$.

Similarly for $i\ge 2$ and $x_i\ne 0$, we also have $x_1= (1-\lambda )x_i$. From this we obtain that $x_i=x_j$.

\noindent  {\bf Case 2: $k$ is even}. If $j\ge 2$ and $x_j\ne 0$, by taking the absolute values of the both sides of Eq. (2) and using the result (1) of this lemma we also have $|x_1|= |(1-\lambda )x_j|$.

Similarly for $i\ge 2$ and $x_i\ne 0$, we also have $|x_1|= |(1-\lambda )x_i|$. From this we obtain that $|x_i|=|x_j|$.
\ep

From Lemma \ref{lem-0} we can obtain the set of all distinct H-eigenvalues and all corresponding H-eigenvectors of the Laplacian tensor $\mathcal {L}$ of the hyperstar $G$ (except for the eigenvalue 1) in the following Proposition \ref{thm-1} (for the case when $k$ is odd) and Proposition \ref{thm-2} (for the case when $k$ is even). The set of all eigenvectors corresponding to the eigenvalue 1 will be given in Proposition \ref{thm-3}.

\begin{Proposition}\label{thm-1}
Let $G=(V,E)$ be a $k$-uniform hyperstar with odd $k\ge 3$ and the size $d\ge 2$, where $V=[n]$, $E=\{e_1,\cdots,e_d\}$, and $d_1=d$ (i.e., the vertex 1 is the heart). Let $\mathcal {L}= \mathcal {D}-\mathcal {A}$ be the Laplacian tensor of $G$. Let
$$f_r(\lambda)=(\lambda -d)(1-\lambda)^{k-1}+r \qquad \qquad  \mbox { for} \quad  r=0,1,\cdots,d$$
Then we have:
\begin{itemize}
\item [(1)] $\lambda \ne 1$ is an H-eigenvalue of $\mathcal {L}$ if and only if it is a real root of the polynomial $f_r(\lambda)$ for some $r\in \{0,1,\cdots,d\}$.
\item [(2)] If $\lambda \ne 1$  is a real root of the polynomial $f_r(\lambda)$, then we can construct all the H-eigenvectors of $\mathcal {L}$ corresponding to $\lambda$ (up to a constant multiple) by going through the following procedure:

\vskip 0.18cm
\noindent {\bf Step 1:} Take $x_1=1-\lambda$.

\vskip 0.18cm
\noindent {\bf Step 2:} Choose any $r$ edges of $G$, take the $x$-values of all the pendant vertices of these $r$ edges to be 1.

\vskip 0.18cm
\noindent {\bf Step 3:} Take the $x$-values of all the other vertices of $G$ to be zero.
\end{itemize}
\end{Proposition}

\noindent {\bf Proof}: (i) Necessity. Let $(\lambda, \mathbf x)$ is an H-eigenpair of $\mathcal {L}$ with $\lambda \ne 1$.

According to the results of Lemma \ref{lem-0}, we call an edge $e$ as $x$-nonzero, if the common $x$-value of all the pendant vertices of $e$ is nonzero. Otherwise this edge is called $x$-zero.

Let $r$ be the number of $x$-nonzero edges of $G$. Then we have $0\le r \le d$. If $r=0$, then $\mathbf x=(1,0,\cdots,0)^T$ is an H-eigenvector corresponding to the eigenvalue $d$, which is the unique root of $f_0(\lambda)$ other than 1. So in the following, we may assume that $1\le r \le d$.

By result (ii) of Lemma \ref{lem-0}, we may assume that $x_i=1$ for all $i\ge 2$ with $x_i\ne 0$ (up to a constant multiple). In this case, we also have $x_1=1-\lambda$.

Now from \reff{hyst-1} we further have $(\lambda -d)x_1^{k-1}=-r$. Combining this with $x_1=1-\lambda$ we obtain that $(\lambda -d)(1-\lambda)^{k-1}+r=0$, which means that $\lambda$ is a real root of the polynomial $f_r(\lambda)$.

Sufficiency part of (i) follows directly from the constructive procedure of result (ii).

(ii)  It is not difficult to verify that any vector $x$ obtained after going through the steps 1-3 will satisfy the \reff{hyst-1} and \reff{hyst-2}, so it is an H-eigenvector corresponding to the H-eigenvalue $\lambda$.
\ep

Now we consider the case when $k$ is even.

\begin{Proposition}\label{thm-2}
 Let $G=(V,E)$ be a $k$-uniform hyperstar with even $k\ge 4$ and the size $d\ge 2$, where $V=[n]$, $E=\{e_1,\cdots,e_d\}$, and $d_1=d$. Let $\mathcal {L}= \mathcal {D}-\mathcal {A}$ be the Laplacian tensor of $G$. Then we have:
\begin{itemize}
\item [(i)] $\lambda \ne 1$ is an H-eigenvalue of $\mathcal {L}$ if and only if it is a real root of the polynomial $f_r(\lambda)$ for some $r\in \{0,1,\cdots,d\}$.
\item [(ii)] If $\lambda \ne 1$  is a real root of the polynomial $f_r(\lambda)$, then we can construct all the H-eigenvectors of $\mathcal {L}$ corresponding to $\lambda$ (up to a constant multiple) by going through the following procedure:

\vskip 0.18cm
\noindent {\bf Step 1:} Take $x_1=1-\lambda$.

\vskip 0.18cm
\noindent {\bf Step 2:} Choose any $r$ edges of $G$, take the $x$-values of all the pendant vertices of these $r$ edges to be $\pm 1$, where the number of $-1$ value in each edge is even.

\vskip 0.18cm
\noindent {\bf Step 3:} Take the $x$-values of all the other vertices of $G$ to be zero.
\end{itemize}
\end{Proposition}

\noindent {\bf Proof}: (i) Necessity. Let $(\lambda, x)$ be an
H-eigenpair of $\mathcal {L}$ with $\lambda \ne 1$.

Let $r$ be the number of $x$-nonzero edges of $G$. Then we have $0\le r \le d$.

By result (ii) of Lemma \ref{lem-0}, we may assume that $x_i=\pm 1$ for all $i\ge 2$ with $x_i\ne 0$ (up to a constant multiple). In this case, we also have $x_1=\pm (1-\lambda)$ by \reff{hyst-2}. We now consider the following two cases:

\noindent  {\bf Case 1: $x_1=(1-\lambda)$ }.  Then from \reff{hyst-2} we have $ \prod_{s\in e(j)\backslash \{1\}}x_s=1$ for $j\ge 2$ and $x_j\ne 0$. Thus from \reff{hyst-1} we further have $(\lambda -d)x_1^{k-1}=-r$. Combining this with $x_1=1-\lambda$ we obtain that $(\lambda -d)(1-\lambda)^{k-1}+r=0$, which means that $\lambda$ is a real root of the polynomial $f_r(\lambda)$.

\noindent  {\bf Case 2: $x_1=-(1-\lambda)$ }. Then from \reff{hyst-2} we have $ \prod_{s\in e(j)\backslash \{1\}}x_s=-1$ for $j\ge 2$ and $x_j\ne 0$. Thus from \reff{hyst-1} we further have $(\lambda -d)x_1^{k-1}=r$. Combining this with $x_1=-(1-\lambda)$ and the hypothesis that $k$ is even, we also obtain that $(\lambda -d)(1-\lambda)^{k-1}+r=0$, which means that $\lambda$ is a real root of the polynomial $f_r(\lambda)$.

Notice that the eigenvectors $x$ constructed in Case 1 and Case 2  only differ by a multiple $-1$, so we only need to consider Case 1.

Sufficiency part of (i) follows directly from the constructive procedure of result (ii).

(ii)  It is not difficult to verify that any vector $x$ obtained after going through the steps 1-3 will satisfy the \reff{hyst-1} and \reff{hyst-2}, so it is an H-eigenvector corresponding to the H-eigenvalue $\lambda$. \ep

Now we construct all the eigenvectors of $\mathcal {\mathcal L}$ corresponding to the eigenvalue 1.

\begin{Proposition}\label{thm-3}
Let $G=(V,E)$ be a $k$-uniform hyperstar with  $k\ge 3$ and the size $d\ge 2$, where $V=[n]$, $E=\{e_1,\cdots,e_d\}$, and $d_1=d$. Let $\mathcal {L}= \mathcal {D}-\mathcal {A}$ be the Laplacian tensor of $G$. Then a nonzero vector $x$ is an eigenvector corresponding to the eigenvalue 1 if and only if $x_1=0$ and the $x$-values of all the pendant vertices of $G$ satisfy the following relation:
\begin{eqnarray}\label{hyst-3}
\sum_{i=1}^d\left ( \prod_{s\in e_i\backslash \{1\}}x_s \right )=0.
\end{eqnarray}
\end{Proposition}

\noindent {\bf Proof}: When $\lambda =1$, \reff{hyst-2} becomes
\begin{eqnarray}\label{hyst-4}
\prod_{s\in e(j)\backslash \{j\}}x_s=0 \qquad \qquad  (j\ge 2).
\end{eqnarray}

Necessity. Suppose that  $\mathbf x$ is an eigenvector corresponding to the eigenvalue 1. If $x_1\ne 0$, then from \reff{hyst-4} we see that each edge of $G$ contains at least two pendant vertices whose $x$-values are zero. From this and the \reff{hyst-1}, we would have $d=1$, a contradiction. So we have that $x_1=0$.

Now $x_1=0$ means that \reff{hyst-1} becomes \reff{hyst-3}. This proves the necessity part.

Sufficiency. It is easy to verify that if $x_1=0$ and the $x$-values of all the pendant vertices of $G$ satisfy the relation \reff{hyst-3}, then $x$ satisfies \reff{hyst-1} and \reff{hyst-2} for $\lambda =1$. Thus $x$ is an eigenvector corresponding to the eigenvalue 1.
\ep

\subsection{Hyperpaths}\label{sec-hp}

In this subsection, we consider a hyperpath of length being $3$ when $k$ is odd.

The next lemma follows from \cite[Theorem 3]{q12a}.
\begin{Lemma}\label{lem-000}
Let $G=(V,E)$ be a $k$-uniform hypergraph and $\mathcal L$ be its Laplacian tensor. If $\lambda\in\mathbb R$ is an H-eigenvalue of $\mathcal L$, then $\lambda\geq 0$.
\end{Lemma}

\begin{Proposition}\label{prop-sph-1}
Let $k$ be odd and $G=(V,E)$ be a $k$-uniform hyperpath with length being $3$. Let $\mathcal L$ be its Laplacian tensor. Then $\lambda\neq 1$ is an H-eigenvalue of $\mathcal L$ if and only if one of the following four cases happens:
\begin{itemize}
\item [(i)] $\lambda=2$ or $\lambda=0$,
\item [(ii)] $\lambda$ is the unique root of the equation $(\lambda-2)(1-\lambda)^{k-1}+1=0$, which is in $(0,1)$,
\item [(iii)] $\lambda$ is the unique root of the equation $(\lambda-2)^2(1-\lambda)^{k-2}-1=0$, which is in $(0,1)$, and
\item [(iv)] $\lambda$ is a real root of the equation $(\lambda-2)^2(1-\lambda)^{k-1}+2\lambda-3=0$ in $(0,2)$.
\end{itemize}
\end{Proposition}

\noindent {\bf Proof.} Suppose that $E=\{\{1,\ldots,k\},\{k,\ldots,2k-1\},\{2k-1,\ldots, 3k-2\}\}$, and $\mathbf x\in\mathbb R^n$ be an H-eigenvector of $\mathcal L$ corresponding to $\lambda\neq 1$.

Let $x_k=\alpha$ and $x_{2k-1}=\beta$. By Lemmas \ref{lem-ch-1} and \ref{lem-ph-1}, we have $x_1=\cdots=x_{k-1}=\frac{1}{1-\lambda}\alpha$ if there are nonzero; $x_{k+1}=\cdots=x_{2k-2}=\pm \sqrt{\frac{\alpha\beta}{1-\lambda}}$ if there are nonzero; and $x_{2k}=\cdots=x_{3k-2}=\frac{1}{1-\lambda}\beta$ if there are nonzero.

The proof is divided into two cases, which contain several sub-cases respectively.

{\bf Case 1} We assume that $x_{k+1}=0$.

(I).
If $x_1=0$, then we must have that either $\lambda=2$ or $\alpha=0$, since $(\lambda-2)\alpha^{k-1}=0$. If $\alpha=0$, then we can assume that $\beta=1$. Thus, either $(\lambda-2)=-\left(\frac{1}{1-\lambda}\right)^{k-1}$ whenever $x_{2k}\neq 0$ or $\lambda=2$. Hence, we have that
either $\lambda=2$ or it is a root of the equation $(\lambda-2)(1-\lambda)^{k-1}=-1$. Let $f(\lambda)=(\lambda-2)(1-\lambda)^{k-1}+1$. We see that $f(0)=-1<0$ and $f(1)=1>0$. Moreover, $f$ is a strictly increasing function in $(-\infty, 1)$ and $(2,\infty)$. We have that $f(\lambda)>0$ in  $(2,\infty)$, since $f(2)=1>0$.
Obviously, $f=0$ does not have a root in $[1,2]$. Thus, it has a unique root, which is in the interval $(0,1)$.

(II).
If $x_1\neq0$, then we have
\begin{eqnarray*}
(\lambda-2)\alpha^{k-1}=-\left(\frac{1}{1-\lambda}\alpha\right)^{k-1}.
\end{eqnarray*}
Since $\alpha\neq 0$ in this case by Lemma \ref{lem-ph-1}, $\lambda$ should be the unique root of the equation $(\lambda-2)(1-\lambda)^{k-1}+1=0$.

The discussion for the cases (i) $x_{2k}=0$, and
(ii) $x_{2k}\neq0$ are similar, and either
$\lambda=2$ or it is the unique root of the equation
$(\lambda-2)(1-\lambda)^{k-1}+1=0$.

{\bf Case 2} We assume that $x_{k+1}\neq 0$.

(I).
If $x_1=0$ and $x_{2k}=0$, then we have
\begin{eqnarray}\label{ph-1}
(\lambda-2)\alpha^{k-1}=-\left(\pm\sqrt{\frac{\alpha\beta}{1-\lambda}}\right)^{k-2}\beta,\;\mbox{and}\;
(\lambda-2)\beta^{k-1}=-\left(\pm\sqrt{\frac{\alpha\beta}{1-\lambda}}\right)^{k-2}\alpha.
\end{eqnarray}
Multiplying the first equality by $\alpha$ and the second by $\beta$, we get that
\begin{eqnarray}\label{ph-3}
(\lambda-2)(\alpha^k-\beta^k)=0.
\end{eqnarray}
If $\lambda>1$, then we have that $\alpha\beta<0$ by Corollary \ref{cor-ph-1}. Thus, the only possibility would be $\lambda=2$ in this case. But $\lambda=2$ contradicts \reff{ph-1}. Hence, in this case we should have that $\lambda<1$. Then, by \reff{ph-3}, we must have $\alpha=\beta\neq 0$ since $k$ is odd. By \reff{ph-1}, we get that
$x_{k+1}$ should be $\sqrt{\frac{\alpha\beta}{1-\lambda}}$ since $k$ is odd and $\lambda<1$. Thus,
$\lambda$ should be a root of the equation $(\lambda-2)^2(1-\lambda)^{k-2}-1=0$ in $(0,1)$. With a similar discussion as that in (I) of {\bf Case 1}. we have that $(\lambda-2)^2(1-\lambda)^{k-2}-1=0$ has a unique root, which is in $(0,1)$.

(II).
If $x_1=0$ and $x_{2k}\neq 0$, then we have
\begin{eqnarray}
(\lambda-2)\alpha^{k-1}&=&-\left(\pm\sqrt{\frac{\alpha\beta}{1-\lambda}}\right)^{k-2}\beta\label{ph-2}\\
(\lambda-2)\beta^{k-1}&=&-\left(\pm\sqrt{\frac{\alpha\beta}{1-\lambda}}\right)^{k-2}\alpha-\left(\frac{1}{1-\lambda}\beta\right)^{k-1}\nonumber.
\end{eqnarray}
Thus, we have
\begin{eqnarray*}
(\lambda-2)\alpha^k=\left[(\lambda-2)+\left(\frac{1}{1-\lambda}\right)^{k-1}\right]\beta^k
\end{eqnarray*}
Let $t:=\frac{\alpha}{\beta}$.
Then
\begin{eqnarray}\label{ph-4}
(\lambda-2)(1-\lambda)^{k-1}t^k=(\lambda-2)(1-\lambda)^{k-1}+1.
\end{eqnarray}
We must have that $\lambda<2$, since $\lambda\leq 2$ by Proposition \ref{prop-ph-3} and $\lambda=2$ cannot be a solution of \reff{ph-4} for any $t\in\mathbb R$.

By squaring the both sides of \reff{ph-2}, we get that
\begin{eqnarray}\label{ph-5}
(\lambda-2)^2(1-\lambda)^{k-2}t^{2k-2}=t^{k-2}.
\end{eqnarray}
By \reff{ph-4} and \reff{ph-5}, $(\lambda,t)$ should be a common solution pair of the polynomial equations
\begin{eqnarray*}
(\lambda-2)(1-\lambda)^{k-1}t^k-(\lambda-2)(1-\lambda)^{k-1}-1&=&0,\\
(\lambda-2)^2(1-\lambda)^{k-2}t^k-1&=&0.
\end{eqnarray*}
Since $k$ is odd, solve $t$ from the second equation, we get that $(\lambda-2)^2(1-\lambda)^{k-1}+2\lambda-3=0$. This, together with Lemma \ref{lem-000}, implies the result (iv).

The discussion for the case $x_{k+1}\neq 0$, $x_1\neq0$ and $x_{2k}=0$ is similar, and the result is the same as the above case.

(III). If $x_1\neq 0$ and $x_{2k}\neq 0$, then we have
\begin{eqnarray}
(\lambda-2)\alpha^{k-1}&=&-\left(\pm\sqrt{\frac{\alpha\beta}{1-\lambda}}\right)^{k-2}\beta-\left(\frac{1}{1-\lambda}\alpha\right)^{k-1}\label{ph-6}\\
(\lambda-2)\beta^{k-1}&=&-\left(\pm\sqrt{\frac{\alpha\beta}{1-\lambda}}\right)^{k-2}\alpha-\left(\frac{1}{1-\lambda}\beta\right)^{k-1}\nonumber.
\end{eqnarray}
Thus, we have
\begin{eqnarray*}
\left[(\lambda-2)+\left(\frac{1}{1-\lambda}\right)^{k-1}\right](\alpha^k-\beta^k)=0.
\end{eqnarray*}
If $\lambda>1$, then we have that $\alpha\beta<0$ by Corollary \ref{cor-ph-1}. Hence, we must have $(\lambda-2)(1-\lambda)^{k-1}+1=0$. But $(\lambda-2)(1-\lambda)^{k-1}+1=0$ has a unique solution in $(0,1)$. Consequently, we must have $\lambda<1$ in this case.

If $\alpha\neq \beta$, then $(\lambda-2)(1-\lambda)^{k-1}+1=0$.
From \reff{ph-6}, we have 
\[
(\lambda-2)(1-\lambda)^{k-1}\alpha^{k-1}+\alpha^{k-1}=-(1-\lambda)^{k-1}\left(\pm\sqrt{\frac{\alpha\beta}{1-\lambda}}\right)^{k-2}\beta.
\]
Consequently, $(1-\lambda)^{k-1}\left(\pm\sqrt{\frac{\alpha\beta}{1-\lambda}}\right)^{k-2}\beta=0$. Hence, $1-\lambda=0$ which is a contradiction to $\lambda<1$. Thus, this case does not happen. 

If $\alpha=\beta$, then by \reff{ph-6} we have that
\begin{eqnarray}\label{ph-7}
\left[(\lambda-2)(1-\lambda)^{k-1}+1\right]^2(1-\lambda)^k=1.
\end{eqnarray}
Note that if $\lambda<0$, then $1-\lambda>1$ and $(\lambda-2)(1-\lambda)^{k-1}+1<-1$, then it cannot be a root of the equation in \reff{ph-7}. If $\lambda\in (0,1)$, then $1-\lambda\in(0,1)$ and $(\lambda-2)(1-\lambda)^{k-1}+1\in (-1, 1)$, then the equation in \reff{ph-7} does not have a root in $(0,1)$. Thus, the unique solution should be $\lambda=0$. \ep

The next lemma says that the degree $d_i$ of the vertex $i$ is a Laplacian H-eigenvalue for all $i\in [n]$.
\begin{Lemma}\label{lem-00}
Let $G=(V,E)$ be a $k$-uniform hypergraph and $\mathcal L$ be its Laplacian tensor. Then $\lambda=d_i$ is an H-eigenvalue of $\mathcal L$ for all $i\in [n]$.
\end{Lemma}

\noindent {\bf Proof.} For any $i\in [n]$, let $\mathbf x\in\mathbb R^{n}$ such that $x_i=1$ and $x_j=0$ for the others. By Definition \ref{def-00}, we have the result. \ep

The next proposition, which follows from Lemma \ref{lem-00}, says that $\lambda=1$ is an H-eigenvalue of the Laplacian tensor of a cored hypergraph and hence a hyperpath.

\begin{Proposition}\label{prop-sph-3}
Let $G=(V,E)$ be a $k$-uniform cored hypergraph and $\mathcal L$ be its Laplacian tensor. Then $\lambda=1$ is an H-eigenvalue of $\mathcal L$.
\end{Proposition}

\begin{Corollary}\label{prop-sph-5}
Let $G=(V,E)$ be a $k$-uniform hyperpath with length being $s\geq 3$ and $\mathcal L$ be its Laplacian tensor. Then $\lambda=1$ is an H-eigenvalue of $\mathcal L$.
\end{Corollary}

\subsection{Hypercycles}\label{sec-hc}

In this subsection, we consider a hypercycle of length being $3$ when $k$ is odd.
\begin{Proposition}\label{prop-sph-2}
Let $k$ be odd and $G=(V,E)$ be a $k$-uniform hypercycle with size being $3$. Let $\mathcal L$ be its Laplacian tensor. Then $\lambda\neq 1$ is an H-eigenvalue of $\mathcal L$ if and only if one of the following four cases happens:
\begin{itemize}
\item [(i)] $\lambda=2$,
\item [(ii)] $\lambda$ is the unique root of the equation $(\lambda-2)^2(1-\lambda)^{k-2}-1=0$, which is in $(0,1)$,
\item [(iii)] $\lambda$ is the unique root of the equation $(\lambda-2)^2(1-\lambda)^{k-2}-2\sqrt[k]{4}=0$, which is in $(0,1)$, and
\item [(iv)] $\lambda$ is a real root of the equation $
\left[(\lambda-2)+\left(\pm\sqrt{1-\lambda}\right)^{k-2}\right](2-\lambda)+2=0$ in [0,1).
\end{itemize}
\end{Proposition}

\noindent {\bf Proof.} Suppose that $E=\{\{1,\ldots,k\},\{k,\ldots,2k-1\},\{2k-1,\ldots, 3k-3,1\}\}$, and $\mathbf x\in\mathbb R^n$ be an H-eigenvector of $\mathcal L$ corresponding to $\lambda\neq 1$.

Let $x_1=\alpha$, $x_{k}=\beta$ and $x_{2k-1}=\gamma$. By Lemmas \ref{lem-ch-1} and \ref{lem-ph-1}, we have $x_2=\cdots=x_{k-1}=\pm \sqrt{\frac{\alpha\beta}{1-\lambda}}$ if there are nonzero; $x_{k+1}=\cdots=x_{2k-2}=\pm \sqrt{\frac{\beta\gamma}{1-\lambda}}$ if there are nonzero; and $x_{2k}=\cdots=x_{3k-3}=\pm \sqrt{\frac{\alpha\gamma}{1-\lambda}}$ if there are nonzero.

(I).
If $x_2\neq 0$, $x_{k+1}=0$ and $x_{2k}=0$, then we have
\begin{eqnarray}\label{sph-1}
(\lambda-2)\alpha^{k-1}=-\left(\pm\sqrt{\frac{\alpha\beta}{1-\lambda}}\right)^{k-2}\beta,\;\mbox{and}\;
(\lambda-2)\beta^{k-1}=-\left(\pm\sqrt{\frac{\alpha\beta}{1-\lambda}}\right)^{k-2}\alpha.
\end{eqnarray}
Multiplying the first equality by $\alpha$ and the second by $\beta$, we get that
\begin{eqnarray}\label{sph-2}
(\lambda-2)(\alpha^k-\beta^k)=0.
\end{eqnarray}
If $\lambda>1$, then we have that $\alpha\beta<0$ by Corollary \ref{cor-ph-1}. Thus, the only possibility would be $\lambda=2$ in this case, since $\lambda\leq 2$. But $\lambda=2$ contradicts \reff{sph-1}. Hence, in this case we should have that $\lambda<1$. Then, by \reff{sph-2}, we must have $\alpha=\beta$ since $k$ is odd. Without loss of generality, we assume that $\alpha=\beta>0$. By \reff{sph-1}, we get that
$x_{k+1}$ should be $\sqrt{\frac{\alpha\beta}{1-\lambda}}$ since $k$ is odd and $\lambda<1$. Thus,
$\lambda$ should be a root of the equation $(\lambda-2)^2(1-\lambda)^{k-2}-1=0$. With a similar discussion as that in (I) of {\bf Case 1} in Proposition \ref{prop-sph-1}, we have that $(\lambda-2)^2(1-\lambda)^{k-2}-1=0$ has a unique root, which is in $(0,1)$.

(II).
If $x_2\neq 0$, $x_{k+1}\neq 0$ and $x_{2k}=0$, then we have
\begin{eqnarray}
(\lambda-2)\alpha^{k-1}&=&-\left(\pm\sqrt{\frac{\alpha\beta}{1-\lambda}}\right)^{k-2}\beta\label{sph-3}\\
(\lambda-2)\beta^{k-1}&=&-\left(\pm\sqrt{\frac{\alpha\beta}{1-\lambda}}\right)^{k-2}\alpha
-\left(\pm\sqrt{\frac{\beta\gamma}{1-\lambda}}\right)^{k-2}\gamma\nonumber\\
(\lambda-2)\gamma^{k-1}&=&-\left(\pm\sqrt{\frac{\beta\gamma}{1-\lambda}}\right)^{k-2}\beta\nonumber.
\end{eqnarray}
Let $\beta=1$, and $s:=\frac{\alpha}{\beta}$ and $t:=\frac{\gamma}{\beta}$. We have
\begin{eqnarray}
(\lambda-2)s^{k-1}&=&-\left(\sqrt{\frac{s}{1-\lambda}}\right)^{k-2}\label{sph-4}\\
(\lambda-2)&=&-\left(\sqrt{\frac{s}{1-\lambda}}\right)^{k-2}s
-\left(\sqrt{\frac{t}{1-\lambda}}\right)^{k-2}t\nonumber\\
(\lambda-2)t^{k-1}&=&-\left(\sqrt{\frac{t}{1-\lambda}}\right)^{k-2}\label{sph-5}.
\end{eqnarray}
Multiplying the first by $s$ and the last by $t$, we have either $\lambda=2$ or $s^k+t^k=1$. $\lambda=2$ contradicts \reff{sph-3}. If $\lambda>1$, then $s^k+t^k=1$ contradicts to the fact that $s<0$ and $t<0$ by Corollary \ref{cor-ph-1}. Thus, $\lambda<1$, $s=t>0$ by \reff{sph-4} and \reff{sph-5}. Since $s^k+t^k=1$, we have $s=t=\sqrt[k]{\frac{1}{2}}$. By \reff{sph-4}, we have that $\lambda$ should be a root of $(\lambda-2)^2(1-\lambda)^{k-2}-2\sqrt[k]{4}=0$. It can be seen that $(\lambda-2)^2(1-\lambda)^{k-2}-2\sqrt[k]{4}=0$ has a unique root, which is in $(0,1)$.

(III).
If $x_2\neq 0$, $x_{k+1}\neq 0$ and $x_{2k}\neq 0$, then we have
\begin{eqnarray}
(\lambda-2)\alpha^{k-1}&=&-\left(\pm\sqrt{\frac{\alpha\beta}{1-\lambda}}\right)^{k-2}\beta
-\left(\pm\sqrt{\frac{\alpha\gamma}{1-\lambda}}\right)^{k-2}\gamma\nonumber\\
(\lambda-2)\beta^{k-1}&=&-\left(\pm\sqrt{\frac{\alpha\beta}{1-\lambda}}\right)^{k-2}\alpha
-\left(\pm\sqrt{\frac{\beta\gamma}{1-\lambda}}\right)^{k-2}\gamma\nonumber\\
(\lambda-2)\gamma^{k-1}&=&-\left(\pm\sqrt{\frac{\beta\gamma}{1-\lambda}}\right)^{k-2}\beta
-\left(\pm\sqrt{\frac{\alpha\gamma}{1-\lambda}}\right)^{k-2}\alpha\nonumber.
\end{eqnarray}
If $\lambda>1$, then $\alpha\beta<0$, $\beta\gamma<0$ and $\gamma\alpha<0$ by Corollary \ref{cor-ph-1}. This is a contradiction. Hence, $\lambda<1$.

Let $\beta=1$, and $s:=\frac{\alpha}{\beta}$ and $t:=\frac{\gamma}{\beta}$. By Lemma \ref{lem-ph-1}, we have $s>0$ and $t>0$. Without loss of generality, we assume that $s\leq 1$ and $t\leq 1$. We have
\begin{eqnarray}
(\lambda-2)s^{k-1}&=&-\left(\pm\sqrt{\frac{s}{1-\lambda}}\right)^{k-2}-\left(\pm\sqrt{\frac{st}{1-\lambda}}\right)^{k-2}t\label{sph-6}\\
(\lambda-2)&=&-\left(\pm\sqrt{\frac{s}{1-\lambda}}\right)^{k-2}s
-\left(\pm\sqrt{\frac{t}{1-\lambda}}\right)^{k-2}t\label{sph-8}\\
(\lambda-2)t^{k-1}&=&-\left(\pm\sqrt{\frac{t}{1-\lambda}}\right)^{k-2}-\left(\pm\sqrt{\frac{st}{1-\lambda}}\right)^{k-2}s\label{sph-7}.
\end{eqnarray}
Multiplying the first by $s$ and the last by $t$, we have that
\begin{eqnarray}\label{sph-9}
(\lambda-2)s^k+\left(\pm\sqrt{\frac{s}{1-\lambda}}\right)^{k-2}s=(\lambda-2)t^k+\left(\pm\sqrt{\frac{t}{1-\lambda}}\right)^{k-2}t
\end{eqnarray}
This, together with \reff{sph-8}, implies that
\begin{eqnarray*}
(\lambda-2)(s^k-1)=(\lambda-2)t^k+2\left(\pm\sqrt{\frac{t}{1-\lambda}}\right)^{k-2}t.
\end{eqnarray*}
Since $s\leq 1$, $\lambda<1$ and $k$ is odd, we have $x_{k+1}=\sqrt{\frac{t}{1-\lambda}}$. Similarly, we have $x_2=\sqrt{\frac{s}{1-\lambda}}$.
Thus, \reff{sph-9} becomes
\begin{eqnarray*}
(\lambda-2)s^k+\left(\sqrt{\frac{s}{1-\lambda}}\right)^{k-2}s=(\lambda-2)t^k+\left(\sqrt{\frac{t}{1-\lambda}}\right)^{k-2}t.
\end{eqnarray*}
It is equivalent to
\begin{eqnarray*}
(s^{\frac{k}{2}}-t^{\frac{k}{2}})\left[(\lambda-2)(s^{\frac{k}{2}}+t^{\frac{k}{2}})+\left(\sqrt{\frac{1}{1-\lambda}}\right)^{k-2}\right]=0.
\end{eqnarray*}
On the other hand, \reff{sph-8} implies that
\begin{eqnarray}\label{sph-10}
2-\lambda=\left(\sqrt{\frac{1}{1-\lambda}}\right)^{k-2}(s^{\frac{k}{2}}+t^{\frac{k}{2}}).
\end{eqnarray}
Hence,
\begin{eqnarray*}
(s^{\frac{k}{2}}-t^{\frac{k}{2}})\left[1-(s^{\frac{k}{2}}+t^{\frac{k}{2}})\right]\left(\sqrt{\frac{1}{1-\lambda}}\right)^{k-2}=0.
\end{eqnarray*}
Thus, either $s=t$ or $s^{\frac{k}{2}}+t^{\frac{k}{2}}=1$.

If $s=t$, then \reff{sph-6} and \reff{sph-10} imply that $\lambda$ and $s$ should be a solution pair of
\begin{eqnarray*}
(2-\lambda)^2(1-\lambda)^{k-2}-4s^k&=&0,\\
\left[(\lambda-2)+\left(\pm\sqrt{\frac{1}{1-\lambda}}\right)^{k-2}\right]s^k+\left(\sqrt{\frac{1}{1-\lambda}}\right)^{k-2}s^{\frac{k}{2}}&=&0.
\end{eqnarray*}
Since $s>0$, we can solve $s$ from the first equation. Then $\lambda$ should be a real root of the equation $
\left[(\lambda-2)+\left(\pm\sqrt{1-\lambda}\right)^{k-2}\right](2-\lambda)2+2=0$ in [0,1), since $\lambda<1$ and $\lambda\geq 0$ by Lemma \ref{lem-000}.

If $s^{\frac{k}{2}}+t^{\frac{k}{2}}=1$, then $\lambda$ should be a root of $(2-\lambda)^2(1-\lambda)^{k-2}-1=0$, which is (ii). \ep


The next corollary, which is a direct consequence of Proposition \ref{prop-sph-3}, says that $\lambda=1$ is also an H-eigenvalue of the Laplacian tensor of a hypercycle.

\begin{Corollary}\label{prop-sph-4}
Let $G=(V,E)$ be a $k$-uniform hypercycle with size being $s\geq 2$ and $\mathcal L$ be its Laplacian tensor. Then $\lambda=1$ is an H-eigenvalue of $\mathcal L$.
\end{Corollary}

\section{Final Remarks}\label{sec-fr}
\setcounter{Theorem}{0} \setcounter{Proposition}{0}
\setcounter{Corollary}{0} \setcounter{Lemma}{0}
\setcounter{Definition}{0} \setcounter{Remark}{0}
\setcounter{Conjecture}{0}  \setcounter{Example}{0} \hspace{4mm} In
this paper, we studied Laplacian H-eigenvalues of cored hypergraphs,
power hypergraphs, and some of their subclasses, such as hyperstars,
hypercycles, hyperpaths and sunflowers.  As the $k$th power of a
tree graph, we have a $k$-uniform hypertree. In 2003, Stevanovi\'c
presented an upper bound for the largest Laplacian eigenvalue of a
tree in terms of the maximum degree. We wonder if this result can be
generalized to hypertrees or not.

\bibliographystyle{model6-names}

\end{document}